\documentclass[preprint,12pt]{elsarticle}

\usepackage{amssymb} 
\usepackage{amsmath} 
\usepackage{bm}
\usepackage{mathrsfs}
\usepackage{soul}
\usepackage{array}
\usepackage{longtable}
\usepackage{caption}
\usepackage{subcaption}
\usepackage{float}
\usepackage{tikz}
\usepackage{parskip}
\usepackage[a4paper, left=2.5cm, right=2.5cm, top=3cm, bottom=3cm]{geometry}
\usepackage[colorlinks = true, linkcolor = blue, urlcolor = blue, citecolor = blue, anchorcolor = blue]{hyperref}
\usepackage{setspace}
\graphicspath{{./figures/}}

\usepackage{lineno}

\usepackage{cleveref}
\usepackage{tcolorbox}
\newtcolorbox{mybox}[1]{fonttitle=\bfseries,title=#1}

\begin{document}


\begin{frontmatter}

\title{The key to the enhanced performance of slab-like topologically interlocked structures with non-planar blocks}

\author[inst1]{Ioannis Koureas}
\author[inst1]{Mohit Pundir}
\author[inst1]{Shai Feldfogel}

\author[inst1]{David S. Kammer\corref{cor1}}
\cortext[cor1]{Corresponding author}
\ead{dkammer@ethz.ch}

\affiliation[inst1]{organization={Institute for Building Materials}, 
            addressline={ETH Zurich},
            city={Zurich},
            country={Switzerland}}

\begin{abstract}
Topologically interlocked structures are assemblies of interlocking blocks that hold together solely through contact. Such structures have been shown to exhibit high strength, energy dissipation, and crack arrest properties. Recent studies on topologically interlocked structures have shown that both the peak strength and work-to-failure saturate with increasing friction coefficient. However, this saturated structural response is only achievable with nonphysically high values of the friction coefficient. For beam-like topologically interlocked structures, non-planar blocks provide an alternate approach to reach similar structural response with friction properties of commonly used materials. It remains unknown whether non-planar blocks have similar effects for slab-like assemblies, and what the achievable structural properties are. Here, we consider slab-like topologically interlocked structures and show, using numerical simulations, that non-planar blocks with wave-like surfaces  allow for saturated response capacity of the structure with a realistic friction coefficient. We further demonstrate that non-planar morphologies cause a non-linear scaling of the work-to-failure with peak strength and result in significant improvements of the work-to-failure and ultimate deflection -- values that cannot be attained with planar-faced blocks. Finally, we show that the key morphology parameter responsible for the enhanced performance of non-planar blocks with wave-like surfaces is the local angle of inclination at the hinging points of the loaded block. These findings shed new light on topologically interlocked structures with non-planar blocks, allowing for a better understanding of their strengths and energy absorption.
\end{abstract}


\end{frontmatter}

\section{Introduction}
\label{sec:introduction}

\begin{figure}[ht!]
    \begin{center}
    \includegraphics[width=\linewidth]{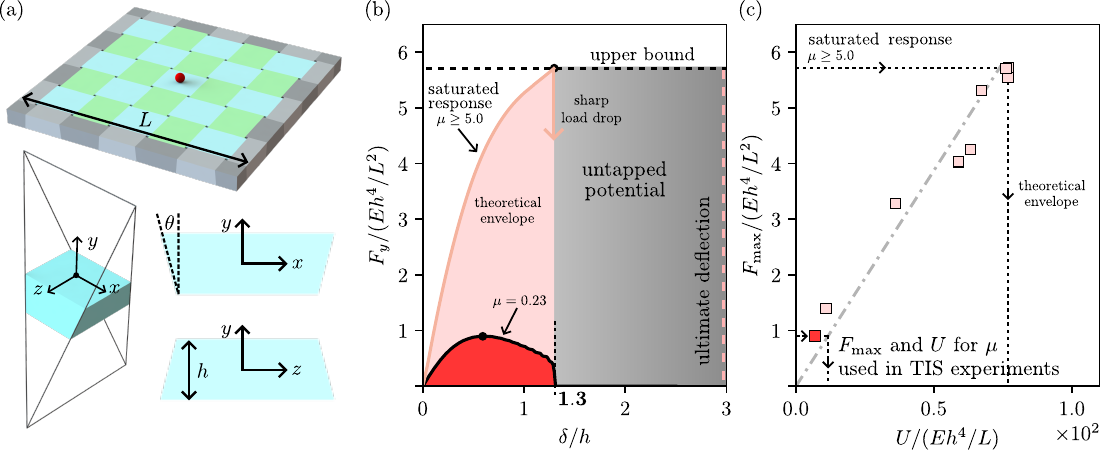}
\caption{(a) Schematic illustration of slab-like TIS with planar-faced blocks along with a single planar block and its cross-sectional area showing the inclination angle $\theta$. (b) Normalized load-deflection curve $F_y / (E \cdot h^4 / L^2)$ vs $\delta/h$, where $E$ is Young's modulus, $h$ the height and $L$ the side length of the panel. The light red curve corresponds to the saturated response of the examined TIS with planar-faced blocks achieved with an (unrealistic) friction coefficient of $\mu = 7.0$. The black curve shows the load-deflection numerically obtained with $\mu=0.23$, a value commonly used for TIS structures in experiments~\cite{Mirkhalaf2019}. The light red shaded region corresponds to the established theoretical envelope and the grey shaded region shows the untapped potential based on \cite{Ullmann2023}. The black points correspond to the peak loads of each curve. (c) Normalized peak load versus work-to-failure. The linear relationship between the two measures shows that for planar blocks, saturation in one measure leads to saturation in the other, achieved with very high values of friction coefficients.}
    \label{fig:fdliterature}
    \end{center}
\end{figure}

Topologically interlocked structures (TIS) are assemblies of interlocking blocks that rely solely on contact along the interfaces to achieve structural integrity~\cite{Dyskin2001, Dyskin2003a}. Slab-like TIS panels (also referred to as plate-like~\cite{Feng2015, RezaeeJavan2017, Williams2021}), a commonly used type of TIS (as shown in \autoref{fig:fdliterature}a), exhibit high strength, energy dissipation~\cite{Dyskin2013}, crack arrest abilities~\cite{Dyskin2003, Dyskin2012, Mirkhalaf2019}, and work-to-failure (the area under the load-deflection curve, elsewhere referred to as toughness~\cite{Mirkhalaf2019, Mirkhalaf2018} or loading energy~\cite{Koureas2022, Koureas2023}).

High work-to-failure $U$ results from the combination of a high peak load $F_{\mathrm{max}}$ and a high ultimate deflection $\delta$. Regarding peak load, it has been shown in \cite{Koureas2022,Feldfogel2022,Zakeri2023} that, in TIS with planar-faced blocks (henceforth referred to as planar blocks), the peak load saturates for high friction coefficients and cannot exceed a well-defined upper bound, see \autoref{fig:fdliterature}b. 
Regarding ultimate deflection, the `envelope' saturated response is commonly capped at relatively small normalized deflection compared to the theoretical upper bound~\cite{Ullmann2023}. For instance, the achieved deflection for a slab-like TIS with 5-by-5 internal blocks is $\delta/h \approx 1.3$ (see \autoref{fig:fdliterature}b), which is $\approx 2.3$ times smaller than the theoretical upper bound $\delta/h=3$.
The limited ability to achieve a high deflection is a direct consequence of the fact that these TIS systems fail in a brittle-like manner as soon as the peak load is reached (see sharp load drop for saturated response in \autoref{fig:fdliterature}b), which results in a linear scaling between the work-to-failure and the peak load (see \autoref{fig:fdliterature}c). As a result, the work-to-failure cannot exceed the value associated to the saturated response (see theoretical envelope in \autoref{fig:fdliterature}b), which points to a large untapped potential for increased work-to-failure (see dark grey region in \autoref{fig:fdliterature}b).

In fact, the untapped potential for work-to-failure is practically even larger because reaching the theoretical envelope requires materials with unrealistically high friction coefficients $\mu$. With realistic $\mu$ in common TIS (\textit{i.e.}, $\mu = 0.23$ as used in ~\cite{Mirkhalaf2019}), the maximally achievable peak load $F_{\mathrm{max}}$ is considerably lower (indicated by the black dot on the black curve in \autoref{fig:fdliterature}b obtained with numerical simulations), and since there is a linear dependence between $F_{\mathrm{max}}$ and $U$, the realistically achievable work-to-failure $U$ is also much lower (see arrows in \autoref{fig:fdliterature}c).



Possible avenues for increasing the work-to-failure in slab-like TIS made from common materials consist of 
1) increasing the peak load to approach the theoretical envelope, and 
2) increasing the ultimate deflection to go beyond this theoretical envelope into the untapped potential.
As we will show here, this is achievable through the use of blocks with non-planar surface morphologies (henceforth referred to as non-planar blocks). Such non-planar morphologies include osteomorphic (double curvature) blocks~\cite{Estrin2003, Dyskin2005, Molotnikov2007, Estrin2011, Tessmann2012, Molotnikov2015, Snijder2016}, double-face osteomorphic blocks~\cite{RezaeeJavan2017,Rezaee2016a, RezaeeJavan2018, RezaeeJavan2020}, and hierarchical morphologies~\cite{Djumas2017}. 

The reason non-planar blocks may outperform planar ones in terms of practically achievable work-to-failure is that they enable a the saturated response with realistic $\mu$, as has been shown for beam-like TIS in \cite{Koureas2023}. If this is also true for slab-like TIS (which are the most common TIS application and the one at the focus of the present study), this would mean an immediate major improvement in the practically achievable work-to-failure by extending the range of accessible values on the linear relation shown in \autoref{fig:fdliterature}c. In addition, if non-planar blocks allowed for less brittle failure, this would correspond to a different scaling between the peak load and the work-to-failure, which, in turn, could lead to further increases in the work-to-failure. Therefore, the degree to which non-planar blocks can help access the untapped work-to-failure potential of slab-like TIS depends on the following questions: Is it possible to reach the saturated response obtained of planar blocks with high $\mu$ by using non-planar morphologies and realistic $\mu$? Is it possible to reach higher ultimate deflections with non-planar blocks than with planar ones? What do such potentially larger deflections mean in terms of achievable work-to-failure and the latter's scaling with the peak load?

While the many studies in the TIS literature that considered non-planar blocks~\cite{RezaeeJavan2017, Estrin2003, Dyskin2005, Molotnikov2007, Estrin2011, Tessmann2012, Molotnikov2015, Snijder2016, Rezaee2016a, RezaeeJavan2018, RezaeeJavan2020, Djumas2017, Xu2020}, provide valuable data and insight, to our knowledge, they did not address the aforementioned questions in the context of slab-like TIS. The objective of the present study is, therefore, to address these questions with a focus on wave-like surface morphologies, and thereby shedding new light on TIS with non-planar blocks, and pointing to ways in which their hitherto untapped work-to-failure potential can be exploited.

\begin{figure}[ht!]
    \begin{center}
    \includegraphics[width=\linewidth]{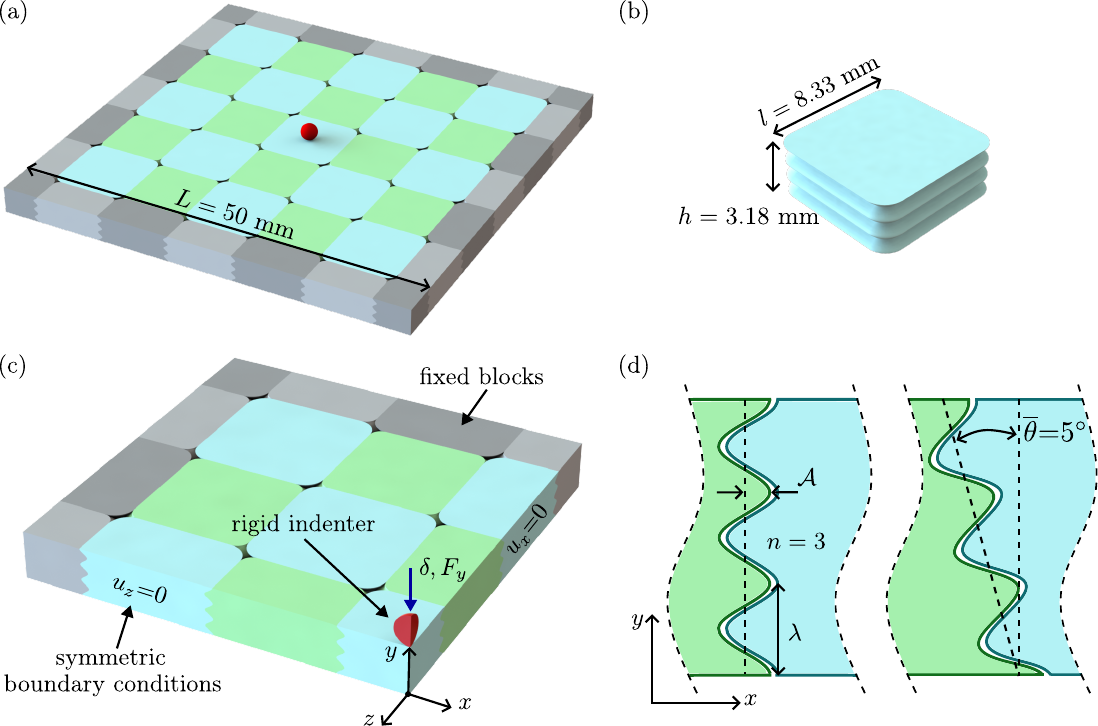}
    \caption{Schematic illustration of (a) a 5-by-5 slab-like topologically interlocked structure with non-planar blocks and structural length $L = 50~\mathrm{mm}$. The immovable peripheral blocks (indicated by grey) support the structure. (b) A single non-planar block with number of oscillations $n=3$, amplitude $\mathcal{A} = 0.25~\mathrm{mm}$, height $h = 3.18~\mathrm{mm}$ and average $l = 8.33~\mathrm{mm}$. (c) Schematic of a quarter of the structure used in the numerical analyses. The rigid indenter is pushed down on the central block by a prescribed displacement $\delta$. (d) Cross-sectional representation of a block's non-planar surface with $\mathcal{A} = 0.25~\mathrm{mm}$, $n = 3$, wavelength $\lambda$, and an average inclination $\overline{\theta} = 0^{\circ}$ (left) and $\overline{\theta} = 5^{\circ}$ (right).}
    \label{fig:tismodel}
    \end{center}
\end{figure}

\section{Problem Statement}
\label{sec:methodology}

\subsection{Physical problem}

In this study, we consider the 5-by-5 slab-like TIS depicted in \autoref{fig:tismodel}a. Each block has length $l = 8.33~\textrm{mm}$ and height $h = 3.18~\textrm{mm}$ (see \autoref{fig:tismodel}b). The immovable peripheral blocks are halves and quarters (in the four corners) and are used to support the structure (see \autoref{fig:tismodel}c). We consider an isotropic and linear elastic material with Young's modulus $E = 18.7~\mathrm{GPa}$, and Poisson's ratio $\nu = 0.2$. A friction coefficient value of $\mu=0.23$ is assumed at the blocks' interface. This configuration was chosen because it allows us to compare the effects of non-planar interfaces with those of planar interfaces~\cite{Mirkhalaf2019, Feldfogel2022} with high friction coefficients~\cite{Feldfogel2022}, as available in the literature.

\subsection{Planar-faced blocks}
\def\subsectionautorefname{Sec.}
\label{sec:planar}

Slab-like TIS with planar-faced blocks as depicted in \autoref{fig:fdliterature}a have been widely used in the literature~\cite{Mirkhalaf2019, Feldfogel2022, Ullmann2023, Khandelwal2012, Khandelwal2014, Khandelwal2015, Mirkhalaf2016}. Here, we use them as a reference for comparison with non-planar block morphologies.

\subsection{Blocks with non-planar morphologies}

Following~\cite{Koureas2023, Djumas2017, Li2011, Li2012, Lin2014, Malik2016, Dalaq2019, Dalaq2020}, we consider sinusoidal wave-like patterns at the interface of the TIS (see \autoref{fig:tismodel}) as prototypes of general non-planar surface morphologies. They are defined by wavelength $\lambda$, number of oscillations $n$, and amplitude $\mathcal{A}$ (see \autoref{fig:tismodel}d) with average inclination angle $\overline{\theta}=0^\circ$ (\textit{i.e.}, simple rectangular block). To understand the effects of each of these variables and to find the governing non-dimensional parameters of non-planar morphologies, we consider different combinations of $n=2, 3, 4$ and $\mathcal{A}=0.025$, $0.05$, $0.1$, $0.15$, $0.2$, $0.25$, $0.3~\mathrm{mm}$.

\subsection{Numerical model}

The 5-by-5 assembly is loaded by applying displacement $\delta$ on a rigid spherical indenter. The peripheral blocks are fully fixed (in all their nodes). Using the structures' symmetry about the $x$- and $z$-axes, we model only a quarter of the TIS (see \autoref{fig:tismodel}c). The sharp edges of the blocks are filleted to reduce stress concentrations. We note that our model does not account for block fracture. Therefore, we define the failure of a structure when the loaded central block loses contact with all its neighbors.

The presented model is solved using the Finite Element software Abaqus~\cite{Systemes2019}. We use a static non-linear solver that accounts for geometrical non-linearities and uses a penalty-based surface-to-surface frictional contact algorithm for handling contact and frictional constraints at interfaces. For describing the frictional forces along an interface and the stick-slip transition, we apply Coulomb's friction law ($T = \mu N$, where $T$ and $N$ are tangential and normal interface tractions, respectively). The blocks are discretized using first-order tetrahedral elements. The element size was determined through a mesh convergence analysis in which the effects of further mesh refinement on the mechanical response were found to be negligible.

\section{Results and Discussion}
\label{sec:results}

\subsection{Load-deflection response for planar-faced blocks}

The load-deflection curves ($F_y-\delta$) with planar blocks, are shown in \autoref{fig:fd}a. Similarly to the Level-Set-Discrete-Element-Method (LS-DEM)~\cite{Feldfogel2022}, our FEM model captures the saturation to the peak load as observed for unrealistically high $\mu$ (i.e., $\mu \ge 5.0$). In addition, for lower values of $\mu$ (i.e., $\mu = 0.6$) the blocks initially stick leading to a peak load and then slip of the central block occurs causing a sharp load drop as shown in \autoref{fig:fd}a. When $\mu$ is decreased even further to experimentally equivalent values~\cite{Mirkhalaf2019}, (i.e., $\mu = 0.23$) the peak load is kept to minimum and failure occurs due to sliding mechanism.

\subsection{The effects of non-planar surface morphologies on load-deflection response}

\begin{figure}[ht!]
    \begin{center}
    \includegraphics[width=\linewidth]{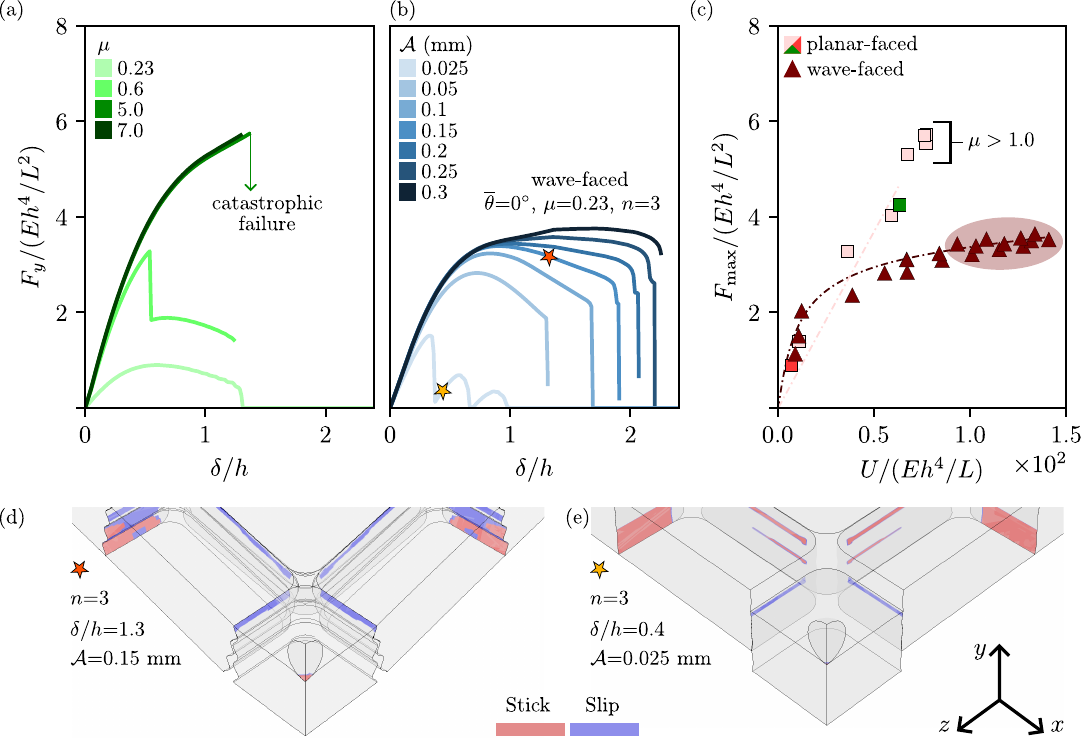}
    \caption{Normalized load $F_y/(E \cdot h^4 / L^2)$ plotted against the normalized prescribed displacement $\delta/h$ for (a) structures with planar-faced blocks obtained with FEM model (solid curves) and those obtained in~\citet{Feldfogel2022} using LS-DEM (dashed curves). (b) Structures with wave-faced blocks and average inclination angle $\overline{\theta} = 0^{\circ}$, oscillation $n = 3$, and different amplitudes $\mathcal{A}$. (c) The normalized peak load $F_{\mathrm{max}}/(E \cdot h^4 / L^2)$ for each curve (for all $n$ and $\mathcal{A}$) is plotted against the normalized work-to-failure $U/(E \cdot h^4 / L)$ for TIS with planar-faced and wave-faced blocks. The dark-shaded cluster denotes the area where $F_{\mathrm{max}}$ saturates for wave-faced blocks. For the complete study of $F_y-\delta$ curves for slab-like TIS with wave-faced blocks examined (\textit{i.e.}, $n=2$ and $n=4$), see \ref{sec:supplementary}. (d, e) Stick-slip mechanism for structures with $n=3$, $\mathcal{A} = 0.15$ mm, and $\mathcal{A} = 0.025$ mm.}
    \label{fig:fd}
    \end{center}
\end{figure}

Turning to TIS with non-planar blocks, we consider the $F_y-\delta$ response with the sinusoidal patterns depicted in \autoref{fig:fd}b. Similarly to planar-faced block structures, we observe a saturation to $F_{\mathrm{max}}$. However, contrary to planar-faced blocks, where unrealistic $\mu$ are needed to reach saturation, with high values of $\mathcal{A}$, the peak load saturates at the relatively low value of $\mu = 0.23$. This means that a theoretical upper bound capacity of slab-like TIS can be obtained through non-planar block morphologies. 

In terms of ultimate deflection and the associated work-to-failure, both are visibly larger with non-planar blocks compared with planar ones. Specifically, for cases with $\mathcal{A} \ge 0.05$~mm, the ultimate deflection approaches $\delta \approx 2.3h$, considerably larger than that with planar-faced blocks ($\delta \approx 1.3h$ as shown in \autoref{fig:fdliterature}b). These findings suggest that non-planar morphologies are inherently superior to planar-faced blocks in terms of (structural) ductility and energy absorption/dissipation.

We note that the results in \autoref{fig:fd} are theoretical since our model does not account for material non-linearity, damage, or fracture. It is likely that in practice, these effects will lead to reduced peak loads and work-to-failure, even for small interlocking angles.

\subsection{Non-linear scaling of load and work-to-failure for TIS with non-planar blocks}

Differently from the case of planar-faced blocks~\cite{Mirkhalaf2019, Feldfogel2022, Mirkhalaf2016} where the work-to-failure $U$ scales linearly with the peak load, see \autoref{fig:fd}c, with non-planar blocks, it scales sub-linearly. This is due to the latter's ability to sustain $F_{\mathrm{max}}$ (or values close to it) for greater deflections (see \autoref{fig:fd}b) compared to planar-faced blocks. In this case, $U$ continues to increase even after $F_{\mathrm{max}}$ reaches saturation (see shaded cluster in \autoref{fig:fd}c). The sub-linear scaling with wave-faced blocks is favorable compared to the linear scaling with planar blocks as it promotes increased work-to-failure. We refer to this ability as pseudo-ductility since it does not involve actual ductile deformation of the bulk material. We attribute it to local enhanced interlocking, that is, the improvement of interlocking between blocks through local non-planar morphologies\footnote{Based on this definition, planar-faced blocks, which lack such local features, have a smaller degree of overall interlocking.}.

To put this result into context with regards to TIS with planar-faced blocks, we highlight two cases. First, we focus on TIS with the same $\mu$ (i.e., 0.23) but with different block types, \textit{i.e.} TIS with planar blocks (highlighted by a red square marker in \autoref{fig:fd}c) with TIS with non-planar blocks (highlighted by shaded cluster in \autoref{fig:fd}c).
TIS with non-planar blocks not only increase significantly the peak load (\textit{i.e.} more than four times), as could have been expected from beam-like observations~\cite{Koureas2023}, but also the work-to-failure is increased by more than a factor of $19$ compared to TIS planar blocks.
This clearly demonstrates the capacity of TIS with non-planar surface morphologies to significantly improve the work-to-failure -- and this while simultaneously increasing the peak load. In addition, if the same peak load needs to be achieved with planar blocks, one would need a friction coefficient as high as 1.0 (highlighted by a green square marker in \autoref{fig:fd}c) and that would still cause a decrease in work-to-failure by more than a factor $2$. Hence, these results suggest that when a material is given for a particular TIS application, introducing non-planar block surface morphology may improve TIS performance both in terms of peak load and work-to-failure.


While TIS with non-planar blocks have been explored previously, they do not exhibit the same improved performance we have identified. This stems from specific surface characteristics, particularly the way the local angle of the surface changes. If we consider osteomorphic blocks as an example, the curvature of the interface promotes smooth sliding along the surface, resulting in sustained work-to-failure, but lower peak loads. In contrast, a wavy pattern, as suggested here with multiple oscillations, introduces curvature that restricts neighboring blocks from sliding extensively across the whole interface. This localized sliding on the curved surface allows us to achieve both high peak loads and work-to-failure and therefore to access the untapped potential.

\subsection{The importance of the local angles of inclination at hinging contact regions}
\label{sec:morphology}

The design of non-planar surface pattern creates enhanced interlocking due to the high local inclination of the interface. While both unrealistically high $\mu$ and non-planar surface patterns can lead to $F_{\mathrm{max}}$ saturation, high $\mu$ alone cannot provide pseudo-ductile behavior and enhanced work-to-failure. This means that, regardless of how high the friction coefficient is, it cannot provide the high levels of work-to-failure provided by high local interlocking angles.

The effect of local enhanced interlocking is significant compared to increasing the angle of inclination ($\theta$) in a planar-faced block because the angle of inclination is limited (design constraint). On the contrary, non-planar blocks with wave-like surfaces can locally introduce very high inclinations, which favor this behavior of promoted work-to-failure while maintaining blocks with attainable shapes. This raises a question regarding the underlying characteristic of the non-planar interfaces that causes the non-linear scaling response.

\begin{figure}[ht!]
    \begin{center}
    \includegraphics[width=\linewidth]{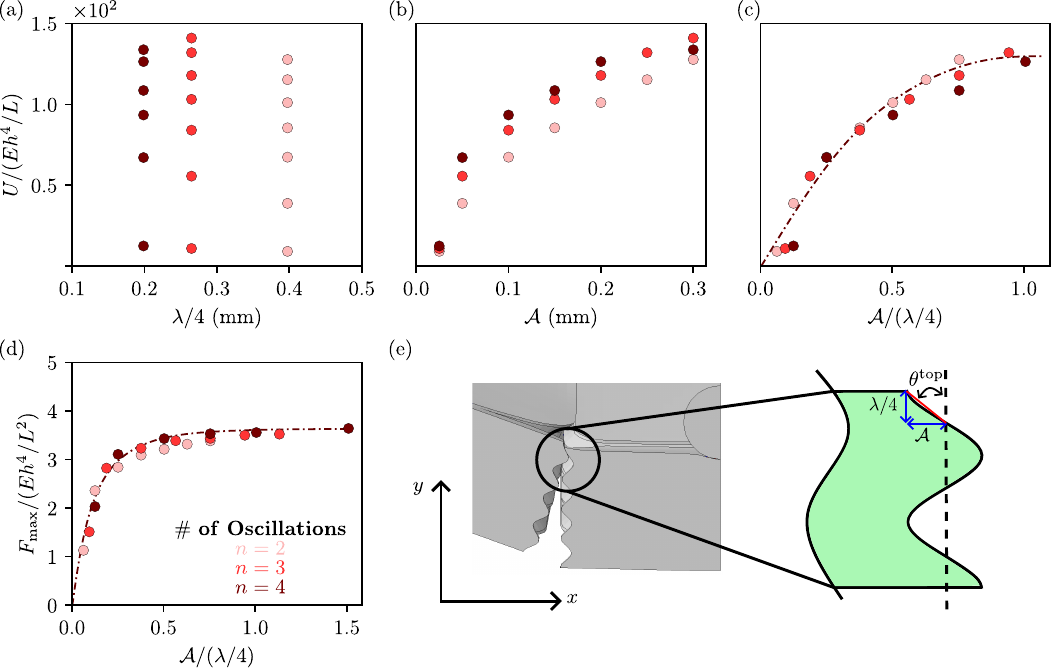}
    \caption{Normalized work-to-failure, plotted against (a) $\lambda/4$ where $\lambda$ is the wavelength of the surface pattern, (b) the amplitude $\mathcal{A}$ and (c) the top inclination $\mathcal{A}/(\lambda/4)$. (d) The normalized peak load for each curve is plotted against the top inclination $\mathcal{A}/(\lambda/4)$. The dashed lines in (c) and (d) have been added as a visual aid to demonstrate the saturation of the work-to-failure and peak load, respectively. (e) Schematic illustration of the top inclination angle between the central and a neighboring block.}
    \label{fig:emorphology}
    \end{center}
\end{figure}

To quantify the effects of each of the parameters governing non-planar block morphologies, we consider $U$ first as a function of $\lambda/4$ where $\lambda$ is the wavelength of the surface pattern \autoref{fig:emorphology}a, and next as a function of waviness amplitude $\mathcal{A}$ in \autoref{fig:emorphology}b. While $\lambda/4$ does not seem to directly influence $U$ as indicated by the large scatter of the data points, there appears to be some correlation with $\mathcal{A}$. Since, we expect the local interlocking to be the main factor, we consider $U$ as a function of $\mathcal{A}/(\lambda/4)$, which represents the average slope at the top of the interface (see \autoref{fig:emorphology}e) as a measure of the degree of interlocking, and observe that the data points collapse and saturate neatly (see \autoref{fig:emorphology}c), as they also do for $F_{\mathrm{max}}$ in \autoref{fig:emorphology}d. This characteristic promotes large rotations and constraints sliding to a specific local area of the overall interface (see \autoref{fig:fd}d). This means that (a) the local angle at the top of the interface of the loaded block (where it is in contact with its neighbors, see $\theta^{\mathrm{top}}$ in \autoref{fig:emorphology}e), not the wavelength or the amplitude independently, is the governing non-dimensional parameter defining the effects of non-planar blocks on the structural response; (b) the local angle is the key to the enhanced work-to-failure and ultimate deflection obtainable with non-planar interfaces; and (c) there is a limit to $\theta^{\mathrm{top}}$ beyond which its further increase does not affect the structural performance. 

We note that high local angles are likely to be associated with high surface tractions and bulk stresses, which can lead to material damage and fracture. This means that the practical applicability of our findings might require technical solutions to address these phenomena.

\subsection{Improving stability for low-amplitude configurations}

\begin{figure}[ht!]
    \begin{center}
    \includegraphics[width=\linewidth]{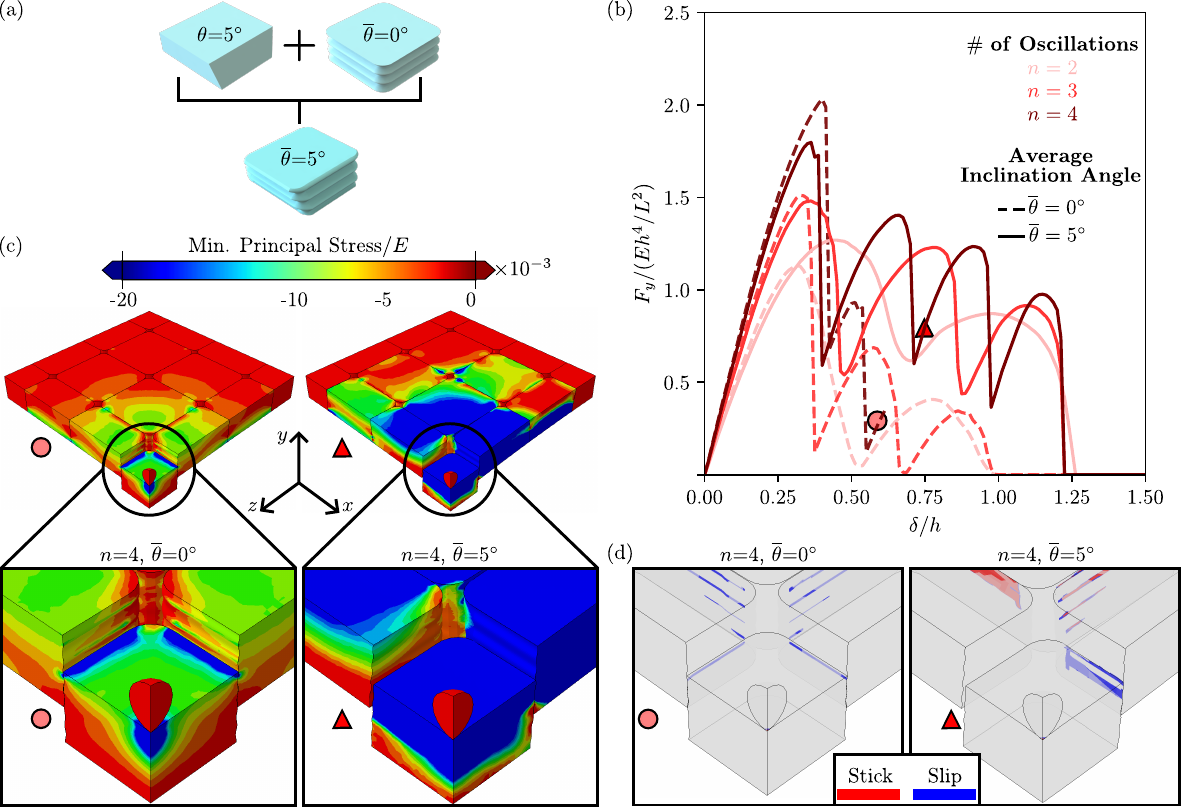}
    \caption{(a) Schematic representation of a wave-faced block with average inclination angle $\overline{\theta} = 5^{\circ}$ as it results from the combination of a planar-faced block with inclination angle $\theta = 5^{\circ}$ and a wave-faced block with average inclination angle $\overline{\theta} = 0^{\circ}$. (b) Normalized load $F_y/(E \cdot h^4 / L)$ against the normalized prescribed displacement $\delta/h$ for structures with different oscillations $n$, and for amplitude $\mathcal{A} = 0.025\mathrm{mm}$. The dashed curves correspond to an average inclination angle $\overline{\theta} = 0^\circ$ and the solid curves to average inclination angle $\overline{\theta} = 5^\circ$. (c) Snapshots of minimum principal stress distribution normalized with Young's modulus $E$ of the structure at specific $\delta/h$ (left) for $n=4$ and $\overline{\theta} = 0^\circ$ and (right) for $n=4$ and $\overline{\theta} = 5^\circ$. (d) Snapshots of the stick-slip behavior at the interfaces of the central block for structures (left) with $\mathcal{A} = 0.025$, $n = 4$, and $\overline{\theta} = 0^\circ$ and (right) $\mathcal{A} = 0.025$, $n = 4$, and $\overline{\theta} = 5^\circ$).}
    \label{fig:truncwavy}
    \end{center}
\end{figure}

Next, we consider the influence of the average inclination of the interface $\overline{\theta}$, see \autoref{fig:truncwavy}a, on the structural response of TIS with wave-faced blocks. Specifically, we focus on how $\overline{\theta} > 0^{\circ}$ can contribute to ameliorating the load-drop instabilities associated with small amplitudes for $\overline{\theta} = 0^{\circ}$ observed in \autoref{fig:fd}b for $\mathcal{A} = 0.025$~mm. Specifically, we observe that the $F_y-\delta$ response undergoes a series of fluctuations characterized by multiple local maxima where each local maximum is successively smaller in magnitude compared to the previous one. This behavior is due to the `pushing-downward' of the central block from the top interface crest to the next one through a sliding mechanism (see \autoref{fig:fd}e). Reconsidering the cases with $\mathcal{A} = 0.025$~mm (see \autoref{fig:fdapp}), but this time with $\overline{\theta} = 5^{\circ}$, we see that, for all the examined number of oscillations $n$, with $\overline{\theta}=5^{\circ}$ the load-drops are smaller, the ultimate deflections increase by roughly $25\%$, and, as a result, the work-to-failure roughly doubles \autoref{fig:truncwavy}b. This speaks to the great contribution of the average angle of inclination to the stability of TIS with non-planar blocks.

To gain a deeper understanding of how the average inclination contributes to the stability of the structure, we consider its effect on the internal force transmission by comparing the distribution of the minimum principal stress with $\overline{\theta} = 0^{\circ}$ and $\overline{\theta} = 5^{\circ}$ (\autoref{fig:truncwavy}c) and the contact forces of the loaded block (\autoref{fig:truncwavy}d). With $\overline{\theta} = 0^{\circ}$, the stress and contact force distribution are symmetrical, meaning that the load is transferred equally in the $yz$ and $xy$ planes. In contrast, with $\overline{\theta} = 5^{\circ}$, they are asymmetric, with more load being transmitted in the $yz$ plane compared with the $xy$ plane as reflected by the larger compressive stresses indicated by blue. This bias indicates that larger $\overline{\theta}$ make one direction (in this case $yz$) stiffer than the perpendicular one, and it suggests that having a stiffer load path contributes to the overall resilience of the structure.

\section{Conclusion}
\label{sec:conclusion}

This study examined how, through the use of blocks with non-planar surface morphologies, the hitherto untapped work-to-failure potential of slab-like topologically interlocked structures (TIS) can be better exploited. For this purpose, we considered blocks with wave-like sinusoidal surface morphologies defined by the number of oscillations $n$, the amplitude $\mathcal{A}$, and the average interface inclination angle $\overline{\theta}$. The effects of these parameters on the mechanical response of the structure were investigated and compared to TIS with planar-faced blocks, shedding new light on their influence in achieving improved structural performance. The main conclusions of our study are:

\begin{itemize}
    \item Unlike TIS with planar blocks, in slab-like TIS with non-planar blocks, the structural response parameters saturate with a relatively low and realistic value of the friction coefficient.
    \item With non-planar blocks, the failure mode is more ductile compared to planar ones. This allows to develop much larger work-to-failure, and results in a sub-linear scaling of the work-to-failure as a function of the peak load.
    \item The key non-dimensional parameter governing non-planar blocks with wave-like surface morphologies, leading to the improved work-to-failure and ultimate deflection, is the local angle of inclination at the hinging contact regions of the structure.
\end{itemize}

In summary, we have shown that introducing local enhanced interlocking to the structure, promotes the system's stability, work-to-failure, and ultimate deflection. This study can be further expanded by considering multi-layer TIS~\cite{YazdaniSarvestani2019} with architected surface patterns and structures with curved geometries (\textit{i.e.}, arches)~\cite{Akleman2020, Sheth2020}. Exploring such configurations could help further improve structural stability, allowing us to obtain the peak load and ultimate deflection simultaneously.
\section*{Acknowledgements}
DSK and MP acknowledge support from the Swiss National Science Foundation under the SNSF starting grant (TMSGI2\_211655).

\clearpage

\appendix

\section{Load-deflection for all amplitudes and oscillations examined}
\label{sec:supplementary}

\begin{figure}[ht!]
    \begin{center}
    \includegraphics[width=\linewidth]{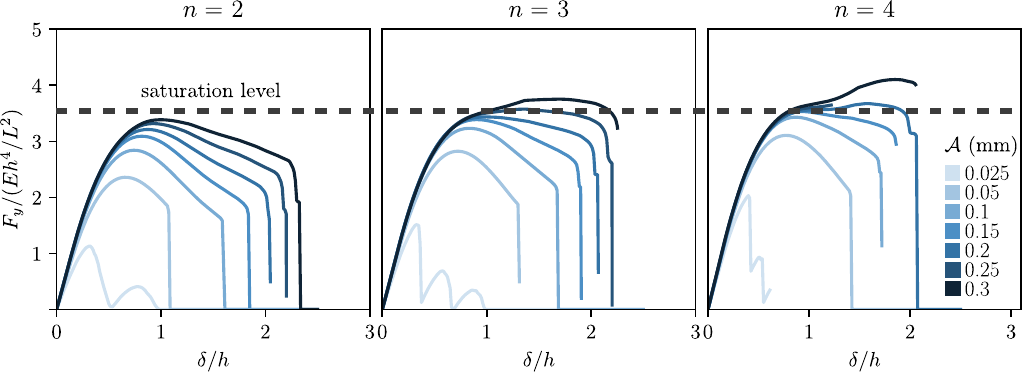}
    \caption{Normalized load-deflection curves for structures with (a) 2, (b) 3 and (c) 4 oscillations $n$, and different amplitudes $\mathcal{A}$. The dashed horizontal line indicates the saturation to the peak load $F_{\mathrm{max}}$.}
    \label{fig:fdapp}
    \end{center}
\end{figure}

\autoref{fig:fdapp} shows that the $F_y - \delta$ response saturates to $F_{\mathrm{max}}$ for all $n$, while further increase of $\mathcal{A}$ promotes pseudo-ductile behavior causing increase in ultimate deflection and work-to-failure. After reaching the saturated response for $F_{\mathrm{max}}$, we have observed significant deformation in the elements, which impacts the accuracy of the $F_y-\delta$ curves (see curves that surpass the saturation level in \autoref{fig:fdapp}b~and~c). This occurrence arises from significantly reducing the number of contact points within the mesh at the interface between the central block and its neighboring ones. The reduced contact points result from high block rotations, resulting in numerical artefacts and recorded load value inaccuracies.

\newpage

\bibliography{cas-refs}

\begin{thebibliography}{41}
\providecommand{\natexlab}[1]{#1}
\providecommand{\url}[1]{\texttt{#1}}
\expandafter\ifx\csname urlstyle\endcsname\relax
  \providecommand{\doi}[1]{doi: #1}\else
  \providecommand{\doi}{doi: \begingroup \urlstyle{rm}\Url}\fi

\bibitem[Dyskin et~al.(2001)Dyskin, Estrin, Kanel-Belov, and
  Pasternak]{Dyskin2001}
A.~V. Dyskin, Y.~Estrin, A.~J. Kanel-Belov, and E.~Pasternak.
\newblock {A new concept in design of materials and structures: Assemblies of
  interlocked tetrahedron-shaped elements}.
\newblock \emph{Scripta Materialia}, 44\penalty0 (12):\penalty0 2689--2694, jun
  2001.
\newblock ISSN 13596462.
\newblock \doi{10.1016/S1359-6462(01)00968-X}.

\bibitem[Dyskin et~al.(2003{\natexlab{a}})Dyskin, Estrin, Kanel-Belov, and
  Pasternak]{Dyskin2003a}
A.~V. Dyskin, Y.~Estrin, A.~J. Kanel-Belov, and E.~Pasternak.
\newblock {Topological interlocking of platonic solids: A way to new materials
  and structures}.
\newblock \emph{Philosophical Magazine Letters}, 83\penalty0 (3):\penalty0
  197--203, 2003{\natexlab{a}}.
\newblock ISSN 09500839.
\newblock \doi{10.1080/0950083031000065226}.

\bibitem[Feng et~al.(2015)Feng, Siegmund, Habtour, and Riddick]{Feng2015}
Yuezhong Feng, Thomas Siegmund, Ed~Habtour, and Jaret Riddick.
\newblock {Impact mechanics of topologically interlocked material assemblies}.
\newblock \emph{International Journal of Impact Engineering}, 75:\penalty0
  140--149, 2015.
\newblock ISSN 0734743X.
\newblock \doi{10.1016/j.ijimpeng.2014.08.003}.
\newblock URL \url{http://dx.doi.org/10.1016/j.ijimpeng.2014.08.003}.

\bibitem[{Rezaee Javan} et~al.(2017){Rezaee Javan}, Seifi, Xu, Ruan, and
  Xie]{RezaeeJavan2017}
Anooshe {Rezaee Javan}, Hamed Seifi, Shanqing Xu, Dong Ruan, and Yi~Min Xie.
\newblock {The impact behaviour of plate-like assemblies made of new
  interlocking bricks: An experimental study}.
\newblock \emph{Materials \& Design}, 134\penalty0 (August):\penalty0 361--373,
  nov 2017.
\newblock ISSN 02641275.
\newblock \doi{10.1016/j.matdes.2017.08.056}.
\newblock URL
  \url{https://linkinghub.elsevier.com/retrieve/pii/S0264127517308146}.

\bibitem[Williams and Siegmund(2021)]{Williams2021}
Andrew Williams and Thomas Siegmund.
\newblock {Mechanics of topologically interlocked material systems under point
  load: Archimedean and Laves tiling}.
\newblock \emph{International Journal of Mechanical Sciences}, 190\penalty0
  (August 2020):\penalty0 106016, 2021.
\newblock ISSN 00207403.
\newblock \doi{10.1016/j.ijmecsci.2020.106016}.
\newblock URL \url{https://doi.org/10.1016/j.ijmecsci.2020.106016}.

\bibitem[Dyskin et~al.(2013)Dyskin, Pasternak, and Estrin]{Dyskin2013}
Arcady Dyskin, Elena Pasternak, and Yuri Estrin.
\newblock {Topological interlocking as a design principle for hybrid
  materials}.
\newblock \emph{8th Pacific Rim International Congress on Advanced Materials
  and Processing 2013, PRICM 8}, 2:\penalty0 1525--1534, 2013.
\newblock \doi{10.1002/9781118792148.ch192}.

\bibitem[Dyskin et~al.(2003{\natexlab{b}})Dyskin, Estrin, Pasternak, Khor, and
  Kanel-Belov]{Dyskin2003}
Arcady~V. Dyskin, Yuri Estrin, Elena Pasternak, Han~C. Khor, and Alexei~J.
  Kanel-Belov.
\newblock {Fracture resistant structures based on topological interlocking with
  non-planar contacts}.
\newblock \emph{Advanced Engineering Materials}, 5\penalty0 (3):\penalty0
  116--119, mar 2003{\natexlab{b}}.
\newblock ISSN 14381656.
\newblock \doi{10.1002/adem.200390016}.
\newblock URL \url{https://onlinelibrary.wiley.com/doi/10.1002/adem.200390016}.

\bibitem[Dyskin et~al.(2012)Dyskin, Pasternak, and Estrin]{Dyskin2012}
Arcady~V. Dyskin, Elena Pasternak, and Yuri Estrin.
\newblock {Mortarless structures based on topological interlocking}.
\newblock \emph{Frontiers of Structural and Civil Engineering}, 6\penalty0
  (2):\penalty0 188--197, 2012.
\newblock ISSN 20952430.
\newblock \doi{10.1007/s11709-012-0156-8}.

\bibitem[Mirkhalaf et~al.(2019)Mirkhalaf, Sunesara, Ashrafi, and
  Barthelat]{Mirkhalaf2019}
Mohammad Mirkhalaf, Amanul Sunesara, Behnam Ashrafi, and Francois Barthelat.
\newblock {Toughness by segmentation: Fabrication, testing and micromechanics
  of architectured ceramic panels for impact applications}.
\newblock \emph{International Journal of Solids and Structures}, 158:\penalty0
  52--65, feb 2019.
\newblock ISSN 00207683.
\newblock \doi{10.1016/j.ijsolstr.2018.08.025}.
\newblock URL \url{https://doi.org/10.1016/j.ijsolstr.2018.08.025
  https://linkinghub.elsevier.com/retrieve/pii/S0020768318303433}.

\bibitem[Mirkhalaf et~al.(2018)Mirkhalaf, Zhou, and Barthelat]{Mirkhalaf2018}
Mohammad Mirkhalaf, Tao Zhou, and Francois Barthelat.
\newblock {Simultaneous improvements of strength and toughness in topologically
  interlocked ceramics}.
\newblock \emph{Proceedings of the National Academy of Sciences of the United
  States of America}, 115\penalty0 (37):\penalty0 9128--9133, 2018.
\newblock ISSN 10916490.
\newblock \doi{10.1073/pnas.1807272115}.

\bibitem[Koureas et~al.(2022)Koureas, Pundir, Feldfogel, and
  Kammer]{Koureas2022}
Ioannis Koureas, Mohit Pundir, Shai Feldfogel, and David~S. Kammer.
\newblock {On the failure of beam-like topologically interlocked structures}.
\newblock \emph{International Journal of Solids and Structures}, 259:\penalty0
  112029, dec 2022.
\newblock ISSN 00207683.
\newblock \doi{10.1016/j.ijsolstr.2022.112029}.
\newblock URL \url{http://arxiv.org/abs/2207.01688
  https://linkinghub.elsevier.com/retrieve/pii/S0020768322004826}.

\bibitem[Koureas et~al.(2023)Koureas, Pundir, Feldfogel, and
  Kammer]{Koureas2023}
Ioannis Koureas, Mohit Pundir, Shai Feldfogel, and David~S. Kammer.
\newblock {Beam-Like Topologically Interlocked Structures With Hierarchical
  Interlocking}.
\newblock \emph{Journal of Applied Mechanics}, 90\penalty0 (8):\penalty0 1--7,
  aug 2023.
\newblock ISSN 0021-8936.
\newblock \doi{10.1115/1.4062348}.
\newblock URL
  \url{https://asmedigitalcollection.asme.org/appliedmechanics/article/90/8/081008/1162998/Beam-Like-Topologically-Interlocked-Structures}.

\bibitem[Feldfogel et~al.(2023)Feldfogel, Karapiperis, Andrade, and
  Kammer]{Feldfogel2022}
Shai Feldfogel, Konstantinos Karapiperis, Jose Andrade, and David~S. Kammer.
\newblock {Scaling, saturation, and upper bounds in the failure of
  topologically interlocked structures}.
\newblock \emph{International Journal of Solids and Structures}, 269\penalty0
  (December 2022):\penalty0 112228, may 2023.
\newblock ISSN 00207683.
\newblock \doi{10.1016/j.ijsolstr.2023.112228}.
\newblock URL \url{http://arxiv.org/abs/2212.11554
  https://linkinghub.elsevier.com/retrieve/pii/S0020768323001257}.

\bibitem[Zakeri et~al.(2023)Zakeri, Safarabadi, and Haghighi-Yazdi]{Zakeri2023}
Milad Zakeri, Majid Safarabadi, and Mojtaba Haghighi-Yazdi.
\newblock {A comprehensive investigation of compressive behavior of
  architectured materials based on topologically interlocking structures:
  Experimental and numerical approaches}.
\newblock \emph{Mechanics Research Communications}, 130\penalty0
  (May):\penalty0 104132, 2023.
\newblock ISSN 00936413.
\newblock \doi{10.1016/j.mechrescom.2023.104132}.
\newblock URL \url{https://doi.org/10.1016/j.mechrescom.2023.104132}.

\bibitem[Ullmann et~al.(2023)Ullmann, Kammer, and Feldfogel]{Ullmann2023}
Silvan Ullmann, David Kammer, and Shai Feldfogel.
\newblock {The deflection limit of slab-like topologically interlocked
  structures}.
\newblock \emph{Journal of Applied Mechanics}, pages 1--14, sep 2023.
\newblock ISSN 0021-8936.
\newblock \doi{10.1115/1.4063345}.
\newblock URL
  \url{https://asmedigitalcollection.asme.org/appliedmechanics/article/doi/10.1115/1.4063345/1166659/The-deflection-limit-of-slab-like-topologically}.

\bibitem[Estrin et~al.(2003)Estrin, Dyskin, Pasternak, Khor, and
  Kanel-Belov]{Estrin2003}
Y.~Estrin, A.~V. Dyskin, E.~Pasternak, H.~C. Khor, and A.~J. Kanel-Belov.
\newblock {Topological interlocking of protective tiles for the space shuttle}.
\newblock \emph{Philosophical Magazine Letters}, 83\penalty0 (6):\penalty0
  351--355, jun 2003.
\newblock ISSN 09500839.
\newblock \doi{10.1080/0950083031000120873}.
\newblock URL
  \url{http://www.tandfonline.com/doi/abs/10.1080/0950083031000120873}.

\bibitem[Dyskin et~al.(2005)Dyskin, Estrin, Pasternak, Khor, and
  Kanel-Belov]{Dyskin2005}
A.~V. Dyskin, Y.~Estrin, E.~Pasternak, H.~C. Khor, and A.~J. Kanel-Belov.
\newblock {The principle of topological interlocking in extraterrestrial
  construction}.
\newblock \emph{Acta Astronautica}, 57\penalty0 (1):\penalty0 10--21, 2005.
\newblock ISSN 00945765.
\newblock \doi{10.1016/j.actaastro.2004.12.005}.

\bibitem[Molotnikov et~al.(2007)Molotnikov, Estrin, Dyskin, Pasternak, and
  Kanel-Belov]{Molotnikov2007}
A.~Molotnikov, Y.~Estrin, A.~V. Dyskin, E.~Pasternak, and A.~J. Kanel-Belov.
\newblock {Percolation mechanism of failure of a planar assembly of interlocked
  osteomorphic elements}.
\newblock \emph{Engineering Fracture Mechanics}, 74\penalty0 (8):\penalty0
  1222--1232, may 2007.
\newblock ISSN 00137944.
\newblock \doi{10.1016/j.engfracmech.2006.07.012}.
\newblock URL
  \url{https://linkinghub.elsevier.com/retrieve/pii/S0013794406002785}.

\bibitem[Estrin et~al.(2011)Estrin, Dyskin, and Pasternak]{Estrin2011}
Y.~Estrin, A.~V. Dyskin, and E.~Pasternak.
\newblock {Topological interlocking as a material design concept}.
\newblock \emph{Materials Science and Engineering C}, 31\penalty0 (6):\penalty0
  1189--1194, 2011.
\newblock ISSN 09284931.
\newblock \doi{10.1016/j.msec.2010.11.011}.
\newblock URL \url{http://dx.doi.org/10.1016/j.msec.2010.11.011}.

\bibitem[Tessmann(2012)]{Tessmann2012}
Oliver Tessmann.
\newblock {Topological Interlocking Assemblies}.
\newblock In \emph{Proceedings of the International Conference on Education and
  Research in Computer Aided Architectural Design in Europe}, volume~2, pages
  211--219, 2012.
\newblock ISBN 9789491207037.
\newblock \doi{10.52842/conf.ecaade.2012.2.211}.
\newblock URL
  \url{http://papers.cumincad.org/cgi-bin/works/paper/ecaade2012_176}.

\bibitem[Molotnikov et~al.(2015)Molotnikov, Gerbrand, Qi, Simon, and
  Estrin]{Molotnikov2015}
A.~Molotnikov, R.~Gerbrand, Y.~Qi, G.~P. Simon, and Y.~Estrin.
\newblock {Design of responsive materials using topologically interlocked
  elements}.
\newblock \emph{Smart Materials and Structures}, 24\penalty0 (2):\penalty0
  025034, feb 2015.
\newblock ISSN 1361665X.
\newblock \doi{10.1088/0964-1726/24/2/025034}.
\newblock URL
  \url{https://iopscience.iop.org/article/10.1088/0964-1726/24/2/025034}.

\bibitem[Snijder et~al.(2016)Snijder, Smits, Bristogianni, and
  Nijsse]{Snijder2016}
A.~Snijder, J.~Smits, T.~Bristogianni, and R.~Nijsse.
\newblock {Design and engineering of a dry assembled glass block pedestrian
  bridge}.
\newblock \emph{Challenging Glass Conference Proceedings - Challenging Glass 5:
  Conference on Architectural and Structural Applications of Glass, CGC 2016},
  \penalty0 (June):\penalty0 547--556, 2016.

\bibitem[{Rezaee Javan} et~al.(2016){Rezaee Javan}, Seifi, Xu, and
  Xie]{Rezaee2016a}
Anooshe {Rezaee Javan}, Hamed Seifi, Shanqing Xu, and Yi~Min Xie.
\newblock {Design of a new type of interlocking brick and evaluation of its
  dynamic performance}.
\newblock \emph{Proceeding of the IASS Annual Symposium 2016}, \penalty0
  (September 2016):\penalty0 1--8, 2016.

\bibitem[{Rezaee Javan} et~al.(2018){Rezaee Javan}, Seifi, Xu, Lin, and
  Xie]{RezaeeJavan2018}
Anooshe {Rezaee Javan}, Hamed Seifi, Shanqing Xu, Xiaoshan Lin, and Yi~Min Xie.
\newblock {Impact behaviour of plate-like assemblies made of new and existing
  interlocking bricks: A comparative study}.
\newblock \emph{International Journal of Impact Engineering}, 116\penalty0
  (February):\penalty0 79--93, jun 2018.
\newblock ISSN 0734743X.
\newblock \doi{10.1016/j.ijimpeng.2018.02.008}.
\newblock URL
  \url{https://linkinghub.elsevier.com/retrieve/pii/S0734743X18300484}.

\bibitem[{Rezaee Javan} et~al.(2020){Rezaee Javan}, Seifi, Lin, and
  Xie]{RezaeeJavan2020}
Anooshe {Rezaee Javan}, Hamed Seifi, Xiaoshan Lin, and Yi~Min Xie.
\newblock {Mechanical behaviour of composite structures made of topologically
  interlocking concrete bricks with soft interfaces}.
\newblock \emph{Materials \& Design}, 186:\penalty0 108347, jan 2020.
\newblock ISSN 02641275.
\newblock \doi{10.1016/j.matdes.2019.108347}.
\newblock URL
  \url{https://linkinghub.elsevier.com/retrieve/pii/S0264127519307853}.

\bibitem[Djumas et~al.(2017)Djumas, Simon, Estrin, and Molotnikov]{Djumas2017}
Lee Djumas, George~P. Simon, Yuri Estrin, and Andrey Molotnikov.
\newblock {Deformation mechanics of non-planar topologically interlocked
  assemblies with structural hierarchy and varying geometry}.
\newblock \emph{Scientific Reports}, 7\penalty0 (1):\penalty0 11844, dec 2017.
\newblock ISSN 20452322.
\newblock \doi{10.1038/s41598-017-12147-3}.
\newblock URL \url{http://www.nature.com/articles/s41598-017-12147-3}.

\bibitem[Xu et~al.(2020)Xu, Lin, and Xie]{Xu2020}
Wenzheng Xu, Xiaoshan Lin, and Yi~Min Xie.
\newblock {A novel non-planar interlocking element for tubular structures}.
\newblock \emph{Tunnelling and Underground Space Technology}, 103\penalty0
  (June):\penalty0 103503, sep 2020.
\newblock ISSN 08867798.
\newblock \doi{10.1016/j.tust.2020.103503}.
\newblock URL \url{https://doi.org/10.1016/j.tust.2020.103503
  https://linkinghub.elsevier.com/retrieve/pii/S0886779820304570}.

\bibitem[Khandelwal et~al.(2012)Khandelwal, Siegmund, Cipra, and
  Bolton]{Khandelwal2012}
S.~Khandelwal, T.~Siegmund, R.~J. Cipra, and J.~S. Bolton.
\newblock {Transverse loading of cellular topologically interlocked materials}.
\newblock \emph{International Journal of Solids and Structures}, 49\penalty0
  (18):\penalty0 2394--2403, sep 2012.
\newblock ISSN 00207683.
\newblock \doi{10.1016/j.ijsolstr.2012.04.035}.
\newblock URL \url{http://dx.doi.org/10.1016/j.ijsolstr.2012.04.035}.

\bibitem[Khandelwal et~al.(2014)Khandelwal, Siegmund, Cipra, and
  Bolton]{Khandelwal2014}
S.~Khandelwal, T.~Siegmund, R.~J. Cipra, and J.~S. Bolton.
\newblock {Scaling of the elastic behavior of two-dimensional topologically
  interlocked materials under transverse loading}.
\newblock \emph{Journal of Applied Mechanics, Transactions ASME}, 81\penalty0
  (3):\penalty0 1--9, 2014.
\newblock ISSN 00218936.
\newblock \doi{10.1115/1.4024907}.

\bibitem[Khandelwal et~al.(2015)Khandelwal, Siegmund, Cipra, and
  Bolton]{Khandelwal2015}
S.~Khandelwal, T.~Siegmund, R.~J. Cipra, and J.~S. Bolton.
\newblock {Adaptive mechanical properties of topologically interlocking
  material systems}.
\newblock \emph{Smart Materials and Structures}, 24\penalty0 (4):\penalty0
  045037, apr 2015.
\newblock ISSN 1361665X.
\newblock \doi{10.1088/0964-1726/24/4/045037}.
\newblock URL
  \url{https://iopscience.iop.org/article/10.1088/0964-1726/24/4/045037}.

\bibitem[Mirkhalaf et~al.(2016)Mirkhalaf, Tanguay, and
  Barthelat]{Mirkhalaf2016}
M.~Mirkhalaf, J.~Tanguay, and F.~Barthelat.
\newblock {Carving 3D architectures within glass: Exploring new strategies to
  transform the mechanics and performance of materials}.
\newblock \emph{Extreme Mechanics Letters}, 7:\penalty0 104--113, jun 2016.
\newblock ISSN 23524316.
\newblock \doi{10.1016/j.eml.2016.02.016}.
\newblock URL
  \url{https://linkinghub.elsevier.com/retrieve/pii/S2352431616300402}.

\bibitem[Li et~al.(2011)Li, Ortiz, and Boyce]{Li2011}
Yaning Li, Christine Ortiz, and Mary~C. Boyce.
\newblock {Stiffness and strength of suture joints in nature}.
\newblock \emph{Physical Review E - Statistical, Nonlinear, and Soft Matter
  Physics}, 84\penalty0 (6):\penalty0 062904, dec 2011.
\newblock ISSN 15393755.
\newblock \doi{10.1103/PhysRevE.84.062904}.
\newblock URL \url{https://link.aps.org/doi/10.1103/PhysRevE.84.062904}.

\bibitem[Li et~al.(2012)Li, Ortiz, and Boyce]{Li2012}
Yaning Li, Christine Ortiz, and Mary~C. Boyce.
\newblock {Bioinspired, mechanical, deterministic fractal model for
  hierarchical suture joints}.
\newblock \emph{Physical Review E - Statistical, Nonlinear, and Soft Matter
  Physics}, 85\penalty0 (3):\penalty0 031901, mar 2012.
\newblock ISSN 15393755.
\newblock \doi{10.1103/PhysRevE.85.031901}.
\newblock URL \url{https://link.aps.org/doi/10.1103/PhysRevE.85.031901}.

\bibitem[Lin et~al.(2014)Lin, Li, Weaver, Ortiz, and Boyce]{Lin2014}
Erica Lin, Yaning Li, James~C. Weaver, Christine Ortiz, and Mary~C. Boyce.
\newblock {Tunability and enhancement of mechanical behavior with additively
  manufactured bio-inspired hierarchical suture interfaces}.
\newblock \emph{Journal of Materials Research}, 29\penalty0 (17):\penalty0
  1867--1875, sep 2014.
\newblock ISSN 20445326.
\newblock \doi{10.1557/jmr.2014.175}.
\newblock URL \url{http://link.springer.com/10.1557/jmr.2014.175}.

\bibitem[Malik and Barthelat(2016)]{Malik2016}
Idris~A. Malik and Francois Barthelat.
\newblock {Toughening of thin ceramic plates using bioinspired surface
  patterns}.
\newblock \emph{International Journal of Solids and Structures},
  97-98:\penalty0 389--399, oct 2016.
\newblock ISSN 00207683.
\newblock \doi{10.1016/j.ijsolstr.2016.07.010}.
\newblock URL
  \url{https://linkinghub.elsevier.com/retrieve/pii/S0020768316301640}.

\bibitem[Dalaq and Barthelat(2019)]{Dalaq2019}
Ahmed~S. Dalaq and Francois Barthelat.
\newblock {Strength and stability in architectured spine-like segmented
  structures}.
\newblock \emph{International Journal of Solids and Structures}, 171:\penalty0
  146--157, 2019.
\newblock ISSN 00207683.
\newblock \doi{10.1016/j.ijsolstr.2019.04.012}.

\bibitem[Dalaq and Barthelat(2020)]{Dalaq2020}
Ahmed~S. Dalaq and Francois Barthelat.
\newblock {Manipulating the geometry of architectured beams for maximum
  toughness and strength}.
\newblock \emph{Materials and Design}, 194:\penalty0 108889, sep 2020.
\newblock ISSN 18734197.
\newblock \doi{10.1016/j.matdes.2020.108889}.
\newblock URL \url{https://doi.org/10.1016/j.matdes.2020.108889
  https://linkinghub.elsevier.com/retrieve/pii/S0264127520304238}.

\bibitem[Systemes(2019)]{Systemes2019}
Dassault Systemes.
\newblock {Abaqus}, 2019.

\bibitem[{Yazdani Sarvestani} et~al.(2019){Yazdani Sarvestani}, Mirkhalaf,
  Akbarzadeh, Backman, Genest, and Ashrafi]{YazdaniSarvestani2019}
H.~{Yazdani Sarvestani}, M.~Mirkhalaf, A.~H. Akbarzadeh, D.~Backman, M.~Genest,
  and B.~Ashrafi.
\newblock {Multilayered architectured ceramic panels with weak interfaces:
  energy absorption and multi-hit capabilities}.
\newblock \emph{Materials and Design}, 167, 2019.
\newblock ISSN 18734197.
\newblock \doi{10.1016/j.matdes.2019.107627}.
\newblock URL \url{https://doi.org/10.1016/j.matdes.2019.107627}.

\bibitem[Akleman et~al.(2020)Akleman, Krishnamurthy, Fu, Subramanian, Ebert,
  Eng, Starrett, and Panchal]{Akleman2020}
Ergun Akleman, Vinayak~R. Krishnamurthy, Chia~An Fu, Sai~Ganesh Subramanian,
  Matthew Ebert, Matthew Eng, Courtney Starrett, and Haard Panchal.
\newblock {Generalized abeille tiles: Topologically interlocked space-filling
  shapes generated based on fabric symmetries}.
\newblock \emph{Computers and Graphics (Pergamon)}, 89:\penalty0 156--166, jun
  2020.
\newblock ISSN 00978493.
\newblock \doi{10.1016/j.cag.2020.05.016}.
\newblock URL
  \url{https://linkinghub.elsevier.com/retrieve/pii/S0097849320300674}.

\bibitem[Sheth and Fida(2020)]{Sheth2020}
Urvi Sheth and Aysha Fida.
\newblock {Funicular structures using topological assemblies}.
\newblock \emph{RE: Anthropocene, Design in the Age of Humans - Proceedings of
  the 25th International Conference on Computer-Aided Architectural Design
  Research in Asia, CAADRIA 2020}, 1:\penalty0 75--84, 2020.
\newblock \doi{10.52842/conf.caadria.2020.1.075}.

\end{thebibliography}

\end{document}